\documentclass[aos,MSNbibl,dvips]{arximspdf}

\doi{10.1214/13-AOS1200} 
\volume{42}
\issue{3}
\pubyear{2014}
\firstpage{850}
\lastpage{871}

\makeatletter
\newproclaim{remark}{Remark}
\newcommand{\rrvert}{\vert}
\newcommand{\llvert}{\vert}
\newcommand{\eqref}[1]{(\ref{#1})}
\def\Var{\operatorname{Var}}
\def\E{\mathrm{E}}
\def\pr{\mathrm{P}}
\def\Cov{\operatorname{Cov}}
\newtheorem{lemma}{Lemma}
\newtheorem{proposition}{Proposition}
\newtheorem{corollary}{Corollary}
\def\S{{\sigma}}
\makeatother

\begin{document}
\begin{frontmatter}

\title{Sharp bounds on the variance in randomized experiments}
\runtitle{Sharp bounds on the variance in randomized experiments}

\begin{aug}
\author[a]{\fnms{Peter M.} \snm{Aronow}\ead[label=e1]{peter.aronow@yale.edu}},
\author[b]{\fnms{Donald P.} \snm{Green}\ead[label=e2]{dpg2110@columbia.edu}}
\and
\author[c]{\fnms{Donald K.~K.} \snm{Lee}\corref{}\ead[label=e3]{donald.lee@yale.edu}}
\runauthor{P.~M. Aronow, D.~P. Green and D.~K.~K. Lee}
\affiliation{Yale University,
Columbia University and Yale University}
\address[a]{P.~M. Aronow\\
Department of Political Science\\
Yale University\\
New Haven, Connecticut 06520\\
USA\\
\printead{e1}}

\address[b]{D. P. Green\\
Department of Political Science\\
Columbia University\\
New York, New York 10027\\
USA\\
\printead{e2}}

\address[c]{D.~K.~K. Lee\\
School of Management (O.R.)\\
Yale University\\
New Haven, Connecticut 06520\\
USA\\
\printead{e3}}

\end{aug}

\received{\smonth{5} \syear{2013}}
\revised{\smonth{12} \syear{2013}}

%
\begin{abstract}
We propose a consistent estimator of sharp bounds on the variance of
the difference-in-means estimator in completely randomized experiments.
Generalizing Robins [\textit{Stat. Med.}  \textbf{7}
(1988) 773--785], our
results resolve a well-known identification problem in causal inference
posed by Neyman [\textit{Statist. Sci.} \textbf{5} (1990) 465--472.
Reprint of the original 1923 paper].
A practical implication of our results is that the upper
bound estimator facilitates the asymptotically narrowest conservative
Wald-type confidence intervals, with applications in randomized
controlled and clinical trials.
\end{abstract}

%


\end{frontmatter}

\section{Introduction}\label{sec1}
We consider the long-standing problem of estimating the variance of the
difference-in-means estimator as applied to a completely randomized
experiment performed on a random sample of size $n$ selected without
replacement from a population of size $N$ under a nonparametric model
of deterministic potential outcomes. It has been known since Neyman
\cite
{neyman23} that neither unbiased nor consistent variance estimation is
generally possible in this setting, due to the fact that the joint
distribution of the potential outcomes can never be fully recovered
from data.

In this paper, we propose an interval estimator that is consistent for
sharp bounds, defined as the smallest interval containing all values of
the variance that are compatible with the observable information. The
upper bound is never larger than and often smaller than conventional
approximations. Our estimator is also applicable to all possible cases
of $N$ and $n$ ($n=N<\infty$, $n<N<\infty$, and $N=\infty$), thus
providing a unified treatment of the problem. In the case where the
outcomes are dichotomous and $n=N<\infty$, our estimator reproduces
Robins \cite{robins88} results. The case $n < N < \infty$
generalizes the
settings considered by prior researchers. Unbiased variance estimation
is not generally possible when $N<\infty$, but our estimator produces
asymptotically sharp bounds. When the population size $N$ is infinite,
our estimator recovers the standard variance point estimator for mean
differences between independent groups \cite{neyman23}.

A practical implication of our work is that it facilitates confidence
intervals that are often narrower than intervals produced by
conventional methods: our upper bound variance estimator may be used to
construct conservative Wald-type confidence intervals for the average
treatment effect. Asymptotically, these intervals are the narrowest
Wald-type intervals that are assured to have at least the nominal
coverage. We illustrate empirical performance using data from an
randomized controlled trial, discuss extensions and provide R code
implementing our estimator.
An implementation in Stata is also available from the authors.

\section{Setting}

Consider a population $U_N$ consisting of $N \geq4$ units. From $U_N$,
$n$ units are randomly sampled into the experimental sample, and the
remaining $N - n$ units are left unsampled. Of the $n$ sampled units,
$m \geq2$ units are randomly assigned to the treatment condition, and
$n - m \geq2$ units are randomly assigned to the control condition.
Let the indicator variable $X^T_i$ be one if unit $i$ is assigned to
the treatment condition, and let the indicator $X^C_i$ be one if unit
$i$ is assigned to the control condition. If $X^T_i = X^C_i = 0$, then
the unit is unsampled. Since units are sampled without replacement,
$X^T_i + X^C_i \leq1$. Without loss of generality, assume an index
ordering $i=1,\ldots,N$ such that those assigned to treatment come
first, $X^T_1,\ldots,X^T_m=1$, and those assigned to control come
after, $X^C_{m+1},\ldots,X^C_n = 1$, and the remaining $N - n$
unsampled units, if any, come last.

Associated with each unit $i$ are two potential outcomes \cite
{neyman23,rubin78} under control and treatment, respectively: $y_{0i}$
and $y_{1i}$. For each unit $i$, the analyst then observes $y_{0i}$
when $X^C_i = 1$ and $y_{1i}$ when $X^T_i = 1$. Given elements $v_i$,
$w_i$ for $i = 1,\ldots,N$, we define the finite population mean $\mu
_N(v)$, finite population variance $\S_N^2(v)$ and finite population
covariance $\S_N(v,w)$, respectively, as
\begin{eqnarray*}
\mu_N(v) &=& \frac{1}{N}\sum_{i=1}^N
v_i,\qquad \S_N^2(v) = \frac
{1}{N}\sum
_{i=1}^N \bigl\{v_i -
\mu_N(v)\bigr\}^2,
\\
\S_N(v,w) &=& \frac{1}{N}\sum_{i=1}^N
\bigl\{v_i - \mu_N(v)\bigr\} \bigl\{w_i -
\mu _N(w)\bigr\}.
\end{eqnarray*}
The average treatment effect for the population $U_N$ is $\tau_N = \mu
_N(y_1) - \mu_N(y_0)$. The difference-in-means estimator of $\tau_N$ is
%
\begin{equation}
\label{eq:est-diffinmeans} \hat\tau_N = \hat\mu_N(y_1) -
\hat\mu_N(y_0) = \frac{1}{m}\sum
_{i=1}^m y_{1i} - \frac{1}{n-m}\sum
_{i=m+1}^n y_{0i},
\end{equation}
with $\E_X (\hat\tau_N) = \tau_N$, where the expectation operator
$\E
_X$ averages over all ${{N}\choose{n}} {{n}\choose{m}}$ possible
treatment assignments.

Our inferential target is the variance of $\hat\tau_N$. Adapting
Freedman \cite{freedman08b}, Proposition~1, the variance is
%
\begin{equation}
\label{eq:varcalc} V_N = \frac{1}{N-1} \biggl\{ \frac{N-m}{m}
\S_N^2(y_{1}) + \frac{N -
(n-m)}{n-m}
\S_N^2(y_{0}) + 2\S_N(y_1,y_0)
\biggr\}.
\end{equation}
The unknown quantities in this expression are $\S_N^2(y_{1})$, $\S
_N^2(y_{0})$ and $\S_N(y_1,y_0)$.
By Cochran \cite{cochran77}, Theorem~2.4, unbiased estimators of $\S
_N^2(y_{1})$ and $\S_N^2(y_{0})$ are
\begin{eqnarray*}
\hat{\S}_N^{2}(y_{1}) &= &\frac{N-1}{N(m-1)}\sum
_{i=1}^m \bigl\{ y_{1i} - \hat
\mu_N(y_1) \bigr\}^2,
\\
\hat{\S}_N^{2}(y_{0}) & =& \frac
{N-1}{N(n-m-1)}
\sum_{i=m+1}^n \bigl\{y_{0i} - \hat
\mu_N(y_0) \bigr\}^2.
\end{eqnarray*}
Since both potential outcomes $y_{0i}$ and $y_{1i}$ for the same unit
can never be observed simultaneously, consistent estimators do not
generally exist for $\S_N(y_1,y_0)$ or for $V_N$ when the population
size $N$ is finite. However, when the population being sampled from is
infinite ($N=\infty$), Neyman \cite{neyman23} noted that the control and
treatment units are effectively sampled independently from their
respective distributions. Hence, the covariance term vanishes, and
$V_N$ is point identified. To see this, let $N\rightarrow\infty$ while
holding $m$ and $n$ fixed so that $V_N \rightarrow\frac{1}{m} \S
_N^2(y_{1}) + \frac{1}{n-m} \S_N^2(y_{0})$, the sampling variance for
the difference of independent means.

\subsection{Neyman \texorpdfstring{\cite{neyman23}}{[13]} approximations when $n = N$}

When $n = N$, the sampling variance of the difference-in-means
estimator reduces to
%
\begin{equation}
\label{eq:true_variance} V_n = \frac{1}{n-1} \biggl\{ \frac{n-m}{m}
\S_n^2(y_{1}) + \frac
{m}{n-m}
\S_n^2(y_{0}) + 2\S_n(y_1,y_0)
\biggr\}.
\end{equation}

Neyman \cite{neyman23} proposed an estimator of $V_n$ that uses the
inequality $2 \S_n(y_1,\break y_0) \leq2  \{ \S_n^2(y_{1})\S
_n^2(y_{0}) \}^{1/2} \leq\S_n^2(y_{1}) + \S_n^2(y_{0})$, by
application of the Cauchy--Schwarz inequality and the inequality of
arithmetic and geometric means. An upper bound estimate for $V_n$ is
obtained by setting $2\S_n(y_1,y_0) = \S_n^2(y_{1}) + \S_n^2(y_{0})$
and substituting $\hat{\S}_n^2(y_{1})$ and $\hat{\S}_n^2(y_{0})$
for $\S
_n^2(y_{1})$ and $\S_n^2(y_{0})$, respectively:
%
\begin{equation}
\hat{V}_n^{a} = \frac{n}{n-1} \biggl\{
\frac{ \hat{\S}_n^2(y_{1})}{m} + \frac{\hat{\S}_n^2(y_{0})}{n-m} \biggr\}.
\end{equation}
Since $\E_X \{\hat{\S}_n^2(y_{1})  \}=\S_n^2(y_{1})$ and
$\E
_X  \{\hat{\S}_n^2(y_{0}) \}=\S_n^2(y_{0})$, $\hat{V}_n^{a}$
is conservative as its bias is nonnegative:
%
\begin{equation}
\label{eq:bias} \E_X \bigl(\hat{V}_n^{a} -
V_n \bigr)= (n-1)^{-1} \bigl\{ \S_n^2(y_{1})
+ \S _n^2(y_{0}) - 2 \S_n(y_1,y_0)
\bigr\} \geq0.
\end{equation}
The estimate $\hat{V}_n^{a}$ is also produced by common estimators that
presuppose sampling from an infinite superpopulation, including
heteroskedasticity-\break robust variance estimators \cite{lin,samiiaronow}
and the standard variance estimate for mean differences between
independent groups \cite{neyman23}. Furthermore, $\hat{V}_n^{a}$ is
known to be unbiased for $V_n$ when effects are constant, as would hold
when there exist no treatment effects whatsoever \cite{gadburyamstat}.
For these reasons, the estimate $\hat{V}_n^{a}$ is often recommended
for the analysis of experimental data \cite{freedmanpp,gerbergreen}.

Neyman \cite{neyman23} also proposed a method for computing bounds on
$V_n$.\break
Given only knowledge of the second moments $\S_n^2(y_1)$ and $\S
_n^2(y_0)$, the\break  sharpest bound on $\S_n(y_1,y_0)$ is given by the
Cauchy--Schwarz inequality:\break  $- \{ \S_n^2(y_1) \S_n^2(y_0)
\}
^{1/2} \leq\S_n(y_1,y_0) \leq \{ \S_n^2(y_1) \S
_n^2(y_0) \}
^{1/2}$. By substituting $\hat{\S}_n^2(y_{1})$ and $\hat{\S
}_n^2(y_{0})$ for $\S_n^2(y_{1})$ and $\S_n^2(y_{0})$, Neyman's bound
estimator is
%
\begin{equation}
\hat{V}_n^{b\pm} = \frac{1}{n-1} \biggl[
\frac{n-m}{m}\hat\S_n^2(y_{1}) +
\frac{m}{n-m}\hat\S_n^2(y_{0}) \pm2 \bigl
\{ \hat{\S}_n^2(y_1) \hat{\S
}_n^2(y_0) \bigr\}^{1/2} \biggr].
\end{equation}
The plus or minus sign is chosen depending on whether an upper or a
lower bound estimate is desired. Neyman recommended choosing $\hat
{V}_n^{b+}$ as a conservative approximation to the true variance, and
suggested that it is ``necessary'' (page 471) to assume that the upper
bound given by the Cauchy--Schwarz inequality holds.

\section{Sharp bounds on $V_N$ given marginal distributions of outcomes}

Under the setting considered, estimates for the marginal distributions
of $y_1$ and $y_0$ exist and can be used to obtain asymptotically sharp
bounds on $V_N$ given the information available. Let $G_N(y) =
N^{-1}\sum_{i=1}^N I(y_{1i} \leq y)$ and $F_N(y) = N^{-1}\sum_{i=1}^N
I(y_{0i} \leq y)$ be the marginal distribution functions of $y_1$ and
$y_0$, respectively. Define their left-continuous inverses as
$G^{-1}_N(u) = \inf\{y\dvtx G_N(y) \geq u \}$ and $F^{-1}_N(u) = \inf\{y\dvtx F_N(y) \geq u \}$. Define also
%
\begin{eqnarray}
\label{eq:exactcovbound} %
\S_N^H
(y_1,y_0) & = & \int_0^1
G_N^{-1}(u) F_N^{-1}(u) \,du -
\mu_N(y_1) \mu_N(y_0),
\nonumber
\\[-8pt]
\\[-8pt]
\nonumber
\S_N^L (y_1,y_0) & = & \int
_0^1 G_N^{-1}(u)
F_N^{-1}(1-u) \,du - \mu _N(y_1)
\mu_N(y_0).
\end{eqnarray}

\begin{lemma}[(Hoeffding)] \label{lemma:covbound} Given only $G_N$ and
$F_N$ and no other information on the joint distribution of
$(y_1,y_0)$, the bound
\[
\S_N^L(y_1,y_0) \leq
\S_N(y_1,y_0) \leq\S_N^H(y_1,y_0)
\]
is sharp. The upper bound is attained if $y_1$ and $y_0$ are
comonotonic, that is, $(y_1,y_0) \sim\{ G_N^{-1}(U),
F_N^{-1}(U)\}$ for a uniform random variable $U$ on $[0,1]$. The lower
bound is attained if $y_1$ and $y_0$ are countermonotonic, \textit{that
is,} $(y_1,y_0) \sim\{ G_N^{-1}(U), F_N^{-1}(1-U)\}$.
\end{lemma}
Lemma~\ref{lemma:covbound} implies that $[V_N^L, V_N^H]$ is the
sharpest interval bound for $V_N$:
\begin{eqnarray*}
V_N^H & = &
\frac{1}{N-1} \biggl\{ \frac{N-m}{m}\S_N^2(y_{1})
+ \frac{N
- (n-m)}{n-m}\S_N^2(y_{0}) + 2
\S_N^H(y_1,y_0) \biggr\},
\\
V_N^L & = & \frac{1}{N-1} \biggl\{ \frac{N-m}{m}
\S_N^2(y_{1}) + \frac{N
- (n-m)}{n-m}
\S_N^2(y_{0}) + 2\S_N^L(y_1,y_0)
\biggr\}.
\end{eqnarray*}

In practice, we observe neither $G_N$ nor $F_N$, but rather their
estimates $\hat G_N(y) = m^{-1}\sum_{i=1}^N X^T_i I(y_{1i} \leq y)$,
$\hat F_N(y) = (n-m)^{-1}\sum_{i=1}^{N} X^C_i I(y_{0i} \leq y)$ and
left-continuous inverses
%
\begin{eqnarray*}
\hat G^{-1}_N(u)
&=& \inf\bigl\{y\dvtx \hat G_N(y) \geq u \bigr\} =
y_{1(\lceil m
u \rceil)},
\\
\hat F^{-1}_N(u) &=& \inf\bigl\{y\dvtx \hat
F_N(y) \geq u \bigr\} = y_{0(m +\lceil
(n-m)u \rceil)},
\end{eqnarray*}
where $y_{1(1)} \leq\cdots\leq y_{1(m)}$ and $y_{0(m+1)} \leq\cdots
\leq y_{0(n)}$ are the ordered observed outcomes, and $\lceil x \rceil$
denotes the smallest integer greater than or equal to $x$. Substituting
$(\hat G_N,\hat F_N)$ for $(G_N,F_N)$ in (\ref{eq:exactcovbound})
yields an interval estimator $ [ \hat\S_N^L(y_1,y_0), \hat\S
_N^H(y_1, y_0) ]$ for $\S_N(y_1,y_0)$:
%
\begin{eqnarray*}
 \hat\S_N^H
(y_1,y_0) & = & \int_0^1
\hat G_N^{-1}(u) \hat F_N^{-1}(u) \,du
- \hat\mu_N(y_1) \hat\mu_N(y_0),
\\
\hat\S_N^L (y_1,y_0) & = & \int
_0^1 \hat G_N^{-1}(u)
\hat F_N^{-1}(1-u) \,du - \hat\mu_N(y_1)
\hat\mu_N(y_0).
\end{eqnarray*}

Let the $[0,1]$-partition $\mathcal{P}_{m,n} = \{p_0, p_1,\ldots, p_{P}
\}$ be the ordered distinct elements of $\{0,1/m,2/m,\ldots,1\} \cup\{
0,1/ ( n-m  ),2/ ( n-m  ),\ldots,1 \}$. Let
$y_{1[i]} = y_{1(\lceil mp_i \rceil)}$ and $y_{0[i]} = y_{0 \{
m+\lceil(n-m) p_i \rceil \}}$. The inverses $\hat G^{-1}_N$ and
$\hat F^{-1}_N$ are piecewise constant since $\hat G^{-1}_N(u) =
y_{1[i]}$ and $\hat F^{-1}_N(u) = y_{0[i]}$ for $u \in(p_{i-1},p_i]$.
In addition, the symmetry $p_i = 1 - p_{P-i}$ implies that $p_i -
p_{i-1} = p_{P+1-i} - p_{P-i}$. Thus, $ [\hat\S_N^L(y_1,y_0),
\hat
\S_N^H(y_1,y_0) ]$ reduces to
%
\begin{eqnarray}
\label{eq:est-covbound} %
 \hat\S_N^H
(y_1,y_0) & = & \sum_{i=1}^P
(p_i-p_{i-1}) y_{1[i]} y_{0[i]} - \hat
\mu_N(y_1) \hat\mu_N(y_0),
\nonumber
\\[-8pt]
\\[-8pt]
\nonumber
\hat\S_N^L (y_1,y_0) & = & \sum
_{i=1}^P (p_i-p_{i-1})
y_{1[i]} y_{0[P+1-i]} - \hat\mu_N(y_1) \hat
\mu_N(y_0),
\end{eqnarray}
where $\hat\mu_N(y_1)$ and $\hat\mu_N(y_0)$ are as defined in (\ref
{eq:est-diffinmeans}).

Substituting $\hat\S_N^2(y_1)$, $\hat\S_N^2(y_0)$, and (\ref
{eq:est-covbound}) for $\{\S_N^2(y_1), \S_N^2(y_0), \S_N(y_1,y_0)\}$ in
the expressions for $V_N^L$ and $V_N^H$, we obtain the interval
estimator $[\hat V_N^L, \hat V_N^H]$ for $V_N$:
%
\begin{eqnarray}
\label{eq:varbound} %
\hat V_N^H
& = & \frac{1}{N-1} \biggl\{ \frac{N-m}{m}\hat\S_N^2(y_{1})
+ \frac{N - (n-m)}{n-m}\hat\S_N^2(y_{0}) + 2\hat
\S_N^H(y_1,y_0) \biggr\},
\nonumber
\\[-8pt]
\\[-8pt]
\nonumber
\hat V_N^L & = & \frac{1}{N-1} \biggl\{
\frac{N-m}{m}\hat\S_N^2(y_{1}) +
\frac{N - (n-m)}{n-m}\hat\S_N^2(y_{0}) + 2\hat
\S_N^L(y_1,y_0) \biggr\}.
\end{eqnarray}
Since Lemma~\ref{lemma:covbound} applies to the sample populations as
well, it follows that $\hat V_N^H$ is never greater than $\hat
{V}_N^{b+}$, and $\hat V_N^L$ is never smaller than $\hat{V}_N^{b-}$. R
code to implement the estimators $\hat V_N^H $ and $\hat V_N^L$ is
presented in Appendix \ref{rcode}.

It is possible to demonstrate that, when outcomes are dichotomous and
$n=N$, our estimator essentially reproduces the estimator proposed by
Robins \cite{robins88}, equation (3), with a slight difference due to finite
population corrections. See Copas~\cite{copas}, Gadbury, Iyer and
Albert \cite{gadburyjspi}, Heckman, Smith and
Clements \cite
{heckmanetal} and Zhang et al.~\cite{zhang} for additional details on
identification of the joint distribution of potential outcomes when
outcomes are dichotomous.

\section{Asymptotic sharpness of interval estimator}

Let $\{ U_N \}_N$ be a nested sequence of finite populations. The
potential outcomes $y_1$ and $y_0$ of each unit are fixed, and hence
the population grows deterministically. As in Isaki and Fuller~\cite
{isakifuller}, we
do not assume that the sequences of treatment assignments are nested;
instead, each $U_N$ hosts its own random assignment. Let $H_N(\cdot
,\cdot)$ be the joint distribution function of $(y_1,y_0)$ for $U_N$.
Under mild conditions on the scaling of $U_N$, the interval estimator
$[\hat V_N^L, \hat V_N^H]$ converges to sharp bounds on $V_N$.

\begin{proposition}\label{prop:asymptotics}
Suppose the following conditions hold as $N\rightarrow\infty$:
\begin{longlist}[1.]
\item[1.]$(m/N,n/N) \rightarrow(\theta\rho,\theta)$ for $\theta\in
(0,1]$ and $\rho\in(0,1)$; 
\item[2.]$H_N$ converges weakly to a limit distribution $H$ with marginals
$G(y) = H(y,\infty)$ and $F(y) = H(\infty,y)$; 
\item[3.]$G_N(y) \rightarrow G(y)$ at any discontinuity point of $G$, and
$F_N(y) \rightarrow F(y)$ at any discontinuity point of $F$; 
%
\item[4.] The sequences of distributions represented by $\{G_N \}_N$ and $\{
F_N \}_N$ are uniformly square-integrable. That is, as $\beta
\rightarrow\infty$, 
\[
\sup_N \Biggl\{ \frac{1}{N} \sum
_{i:y_{1i}^2\geq\beta}^N y_{1i}^2 \Biggr\},\qquad
\sup_N \Biggl\{ \frac{1}{N} \sum
_{i:y_{0i}^2\geq\beta}^N y_{0i}^2 \Biggr\}
\rightarrow0.
\]
\end{longlist}

Then for the collection $\mathcal{H}$ of all bivariate distributions
with marginals $G$ and $F$, the moments of each $h\in\mathcal{H}$
exist up to second order and
\begin{eqnarray*}
NV_N^H &\rightarrow&\frac{1-\theta\rho}{\theta\rho}
\Var_H(y_1) + \frac
{1 - \theta(1-\rho)}{\theta(1-\rho)} \Var_H(y_0)
+ 2 \sup_{h \in
\mathcal{H}} \Cov_h (y_1,y_0),
\\
NV_N^L &\rightarrow&\frac{1-\theta\rho}{\theta\rho}
\Var_H(y_1) + \frac
{1 - \theta(1-\rho)}{\theta(1-\rho)} \Var_H(y_0)
+ 2 \inf_{h \in
\mathcal{H}} \Cov_h (y_1,y_0).
\end{eqnarray*}
Moreover, $(\hat V_N^H - V_N^H, \hat V_N^L - V_N^L) = o_P(1/N)$.
\end{proposition}
\begin{remark}\label{re1}
Condition 3 is used to
establish the
functional convergence of $(G_N,F_N)$ to $(G,F)$. When the units of
$U_N$ are independent and identically distributed samples from a
superpopulation, the condition holds with probability one because of
the strong law of large numbers. The condition is also satisfied if $G$
and $F$ are continuous, regardless of whether or not the units come
from a superpopulation. We thank Professor A.~W. van der Vaart for
suggesting the latter as an alternate sufficient condition for
convergence, which subsequently inspired condition~3.
\end{remark}

\begin{remark}\label{re2} Given condition 2, any convergence
of the marginal second moments of $H_N$ to those of $H$ (should they
exist) necessarily implies condition 4. Thus, the condition
is the weakest possible complement to conditions 1--3.
\end{remark}

 \begin{remark}\label{re3} If condition 4 of Proposition~\ref
{prop:asymptotics} is strengthened to require that $y_1$ and $y_0$ be
bounded, then higher order rates of convergence can be obtained, namely
that $\pr( N| \hat V_N^H - V_N^H | > \varepsilon$) and $\pr( N| \hat
V_N^L - V_N^L | > \varepsilon$) are both of order $\mathcal{O}(1/N)$.
Interested readers are referred to Proposition~\ref{prop:rate} in the \hyperref[app]{Appendix}.
\end{remark}

\textit{Outline of proof}. The random treatment assignment process can
be expressed as a triangular array $\mathcal{X}$ where the $N$th row
$(\mathcal{X}_{N,1},\ldots,\mathcal{X}_{N,N}) = \{
(X^T_1,X^C_1),\ldots
, (X^T_N,X^C_N) \}$ is the treatment/control assignment for population
$U_N$. Since the treatment/control assignment for $U_{N+1}$ is not
related to that for $U_N$, each row of $\mathcal{X}$ is a random vector
of a different probability space. As a result, the sequence of random
distribution functions $(\hat G_N,\hat F_N)$ do not share a common
probability space. However, by treating $(\hat G_N,\hat F_N)$ as random
elements taking values in the product space of \textit{c\`adl\`ag}
functions $D([-\infty,\infty],\mathbb{R})^2$ endowed with the uniform
metric, we show that $(\hat G_N,\hat F_N) \rightarrow(G,F)$ in
probability. It then follows from the Skorohod representation that
there exists a sequence of random elements $(\hat G_N',\hat F_N')$
defined on a common probability space that has the same law as $(\hat
G_N,\hat F_N)$. Moreover, $(\hat G_N',\hat F_N')$ converges to $(G,F)$
almost everywhere. Pathwise convergence of the moments of $(\hat
G_N',\hat F_N')$ then implies probabilistic convergence of the moments
of $(\hat G_N,\hat F_N)$ to the desired result. We refer the reader to
the \hyperref[app]{Appendix} for details of the formal argument.

\section{Confidence intervals for \texorpdfstring{$\tau_N$}{tau N}}

The upper bound estimator $\hat{V}_N^H$ may be used as a basis for
Wald-type confidence intervals for the average treatment effect. The
proof of the following corollary follows directly from Freedman
\cite{freedman08b}, Theorem~1, and associated remarks.
%
\begin{corollary}
Suppose that the support of $H$ is nonsingular and that conditions 1--3 of Proposition~\ref
{prop:asymptotics} hold. Suppose in\vadjust{\goodbreak} addition that condition 4 is strengthened to require uniformly bounded third moments:
\[
\sup_N \Biggl\{ \frac{1}{N}\sum
_{i=1}^N | y_{1i} |^3 \Biggr
\},\qquad \sup_N \Biggl\{ \frac{1}{N}\sum
_{i=1}^N | y_{0i} |^3 \Biggr
\} < \infty.
\]
Then
\[
\frac{\hat{\tau}_N - \tau_N}{ (\gamma\hat V_N^H  )^{1/2}}
\]
converges weakly to the standard normal distribution where $\gamma=
\lim_N (NV_N)/\break  \lim_N (N V_N^H) \leq1$.
\end{corollary}
\begin{remark}\label{re4}
As $\hat V_N^H$ is consistent for the sharp upper bound
on $V_N$, then given large $N$, a confidence interval constructed as
$\hat\tau_N \pm z_{1-\alpha/2} (\hat V_N^H)^{1/2}$ is asymptotically
the narrowest Wald-type confidence interval assured to have at least
the nominal coverage.
\end{remark}

\section{Application}

We consider the randomized controlled trial reported by Harrison and
Michelson \cite
{michelson}, which assessed the intention-to-treat effects of an
experimental phone call on donations to a nonprofit gay rights
organization. The control phone call script contained a standard
appeal. The experimental phone call script included an additional
sentence that revealed the sexual orientation of the volunteer caller.
The finite population $U_N$, which was not selected from any broader
population, contains $N = n = 1561$ subjects, $m = 781$ of whom were
randomly assigned to receive the experimental phone call. Outcomes were
measured in terms of US dollars (USD) received per subject, ranging
from $\$0$ to $\$150$. The mean donation given by subjects assigned to
control was $\hat\mu_N(y_0) = \$1.397$, and the mean donation given by
subjects assigned to treatment was $\hat\mu_N(y_1) = \$0.715$, yielding
the difference-in-means estimate $\hat\tau_N = -\$0.682$.

\begin{table}[b]
\caption{Variance estimates and confidence intervals for Harrison and
Michelson \protect\cite{michelson}}\label{tab1}
\begin{tabular*}{\textwidth}{@{\extracolsep{\fill}}lcc@{}}
\hline
& \textbf{Variance} & \multicolumn{1}{c@{}}{\textbf{95\% confidence}} \\
& \textbf{estimate (USD$^2$)} & \multicolumn{1}{c@{}}{\textbf{interval for} $\bolds{\tau_N}$} \\
\hline
Conventional ($\hat{V}_N^{a}$) & 0.199 & $(-\$1.555, \$0.192)$ \\
Neyman upper bound ($\hat{V}_N^{b+}$) & 0.196 & $(-\$1.548, \$0.185)$ \\[2pt]
Neyman lower bound ($\hat{V}_N^{b-}$) & 0.003 & N$/$A \\[2pt]
Sharp upper bound ($\widehat V_N^H$) & 0.186 & $(-\$1.528, \$0.165$) \\[2pt]
Sharp lower bound ($\widehat V_N^L$) & 0.098 & N$/$A \\
\hline
\end{tabular*}
\end{table}

In Table~\ref{tab1}, we report the variance estimates and confidence intervals
associated with Neyman's\vadjust{\goodbreak} approximations and our proposed estimator. We
find, as expected, that our estimates are sharper than Neyman's
approximations. Compared to the conventional variance estimator $\hat
{V}_N^{a}$, we find that our upper bound estimator yields a 7\%
reduction in the nominal variance. Importantly, if using $\hat
{V}_N^{a}$ as a basis for conservative inference, one would need over
100 additional subjects in order to achieve the same nominal variance
as that of our proposed upper bound estimate $\hat V_N^H$, all else
equal. Similarly, if using $\hat{V}_N^{b+}$, one would need over 75
additional subjects to achieve the nominal variance of $\hat V_N^H$.

\subsection{Simulations}

We use the data from \cite{michelson} to assess the operating
characteristics of the upper bound estimators and associated Wald-type
confidence intervals. These characteristics depend on the underlying
joint distribution of potential outcomes, which cannot be directly
observed and are instead hypothesized as part of these simulations. We
thus impute the missing potential outcomes (potential control outcomes
for treatment subjects, and potential treatment outcomes for control
subjects) by asserting varying hypotheses about treatment effects. We
simulate 25 million random assignments and, for each of these random
assignments, compute the upper bound variance estimates $\hat
{V}_N^{a}$, $\hat{V}_N^{b+}$ and $\hat V_N^H$, and associated
confidence intervals that would have been obtained. For the collection
of 25 million simulations, we calculate the mean variance estimate, the
mean width of the associated 95\% confidence intervals for $\tau_N$ and
the fraction of simulated confidence intervals covering $\tau_N$.

The first hypothesis that we evaluate is the sharp null hypothesis of
no effect whatsoever. This hypothesis, denoted ``Sharp Null,'' assumes
that $y_{0i} = y_{1i}$ for all~$i$. Under the Sharp Null, the treatment
effect estimator variance is $0.199$ USD$^2$. As can be seen in Table~\ref{tab2}, Neyman's estimators predictably perform well since they implicitly
assume that the outcomes are perfectly correlated: the bias \eqref
{eq:bias} for $\hat{V}_N^{a}$ is zero because $\S_N^2(y_1) = \S
_N^2(y_0) = \S_N(y_1,y_0)$. Due to the nonlinearity of the square root
function, the Cauchy--Schwarz inequality implies that $\hat{V}_N^{b+}$
has nonpositive bias ($-0.007$ USD$^2$). The 95\% confidence intervals
associated with $\hat{V}_N^{a}$ and $\hat{V}_N^{b+}$ have coverage of
95.2\% and 94.1\%, respectively (the former is not exactly 95\% because
the sampling distribution of $\tau_N$ is not perfectly normal). Because
$\hat V_N^H \leq\hat V_N^{b+}$, $\hat V_N^H$ is slightly more
negatively biased ($-0.010$ USD$^2$) and has lower coverage ($93.7\%$)
than $\hat V_N^{b+}$.

\begin{table}
\caption{Simulated variance estimator properties under varying
treatment effect hypotheses for Harrison and
Michelson \protect\cite{michelson}, using 25 million
simulated random assignments each}\label{tab2}
\begin{tabular*}{\textwidth}{@{\extracolsep{\fill}}lcccc@{}}
\hline
&  & \textbf{Mean var.} &
\textbf{Mean 95\%} &
\textbf{Coverage} \\
\textbf{Effect hypothesis} &\textbf{Variance estimator}& \textbf{estimate} & \textbf{CI width} & \multicolumn{1}{c@{}}{\textbf{for} $\bolds{\tau_N}$} \\
\hline
Sharp Null & Conventional ($\hat{V}_N^{a}$) & 0.199 & 1.747 & 95.2\% \\
(True Var.: 0.199) &Neyman upper bound ($\hat{V}_N^{b+}$) & 0.193 &
1.724 & 94.1\% \\[2pt]
& Sharp upper bound ($\hat V_N^H$) & 0.189 & 1.703 & 93.7\% \\[2pt]
Heterogeneity A & Conventional ($\hat{V}_N^{a}$) & 0.279 & 2.067 &
96.7\% \\[2pt]
(True Var.: 0.238) &Neyman upper bound ($\hat{V}_N^{b+}$) & 0.268 &
2.028 & 95.9\% \\[2pt]
& Sharp upper bound ($\hat V_N^H$) &0.258 & 1.987 & 95.4\% \\[2pt]
Heterogeneity B & Conventional ($\hat{V}_N^{a}$) & 0.244 & 1.933 &
97.4\% \\[2pt]
(True Var.: 0.186) &Neyman upper bound ($\hat{V}_N^{b+}$) &0.226 &
1.860 & 96.5\% \\[2pt]
& Sharp upper bound ($\hat V_N^H$) & 0.214 & 1.809 & 96.0\%\\
\hline
\end{tabular*}
\end{table}

We next consider two hypotheses that embed treatment effect
heterogeneity, denoted ``Heterogeneity A'' and ``Heterogeneity B.'' Under
Heterogeneity~A, we assume that the sharp null hypothesis holds, with
the exception of 10 subjects who had an observed $y_{0i} = 0$ USD under
control. For these 10 subjects, we assume that $y_{1i} = 100$ USD.
Under Heterogeneity A, the treatment effect estimator variance is
$0.238$ USD$^2$\vadjust{\goodbreak} and, as expected, all variance estimators are
conservative (positively biased). However, the bias, confidence
interval widths, and coverage for $\tau_N$ are all improved when $\hat
V_N^H$ is used in place of either of Neyman's estimators. In
formulating the Heterogeneity B hypothesis, we assume that
Heterogeneity A holds, but, in addition, for all 6 subjects under
treatment with an observed $y_{1i} \geq50$ USD, we assume that $y_{0i}
= 0$ USD. Under Heterogeneity B, the treatment effect estimator
variance is $0.186$ USD$^2$ and, again, while all estimators are
conservative, $\hat V_N^H$ improves over Neyman's estimators.

In Appendix \ref{illustrative}, we further explore the relative
performance of the upper bound estimates under varying assumptions
about the distribution of potential outcomes. Using the Beta
distribution family as an example to represent varying shapes of
marginal treatment and control distributions, we show that it is
possible for $\hat V_N^H$ to materially outperform $\hat V_N^a$ and
$\hat V_N^{b+}$ as the two marginals diverge in shape. Our simulations
therefore illustrate how $\hat V_N^H$ can improve upon Neyman's bounds
under effect heterogeneity.

\section{Discussion}

The proposed variance estimator may also be extended to alternative
designs. For block-randomized designs where the number of units per
block grows asymptotically large, Proposition~\ref{prop:asymptotics}
holds within each block, and thus calculation of the overall variance
is straightforward. In cluster-randomized designs with equally-sized
clusters, the proposed estimator may be used with the unit of analysis
being the cluster and the outcome being the cluster mean. It is also
straightforward to adapt the estimator to completely randomized
experiments with multiple treatments, which may be shown to be
logically equivalent to sampling from a broader population. In
addition, we note that our result can be generalized to characterize
estimation error for arbitrary target populations within the sampling
frame (e.g., unsampled units).

Finally, we remark on the scope of our findings, as our results
presuppose deterministic potential outcomes. When the potential
outcomes are stochastic, the total variance is greater than the
conditional variance (\ref{eq:varcalc}) because of the additional
stochasticity. If one sought to estimate the total variance or bounds
thereof, additional structure would need to be imposed on the
stochastic process (e.g., independence across units and finite
variances); otherwise it is possible for the identification set to be unbounded.

\begin{appendix}\label{app}
\section{Proofs}\label{proofs}

\begin{pf*}{Proof of Lemma~\ref{lemma:covbound}}
Let $H_N(y_1,y_0)$
be the joint distribution function of $(y_1,y_0)$, and define two other
distributions $H_N^H(y_1,y_0) = \min\{G_N(y_1),\break F_N(y_0)\}$ and
$H_N^L(y_1,y_0) = \max\{0,G_N(y_1)+F_N(y_0)-1\}$. All three
distributions have the same marginals $G_N$ and $F_N$. Defining $\E_Q$
as the expectation operator with respect to a distribution $Q$, a
result by Hoeffding, recounted in Tchen \cite{tchen1980}, shows that
\[
\E_{H_N^L} (y_1 y_0) \leq\E_{H_N}
(y_1y_0) \leq\E_{H_N^H} (y_1
y_0).
\]
Since $\{G_N^{-1}(U), F_N^{-1}(U)\} \sim H_N^H$ and $\{G_N^{-1}(U),
F_N^{-1}(1-U)\} \sim H_N^L$, the lower and upper bounds are equivalent to
\begin{eqnarray*}
\E_{H_N^H} (y_1 y_0) &=& \int
_0^1 G_N^{-1}(u)
F_N^{-1}(u) \,du,
\\
\E_{H_N^L} (y_1 y_0) &= &\int
_0^1 G_N^{-1}(u)
F_N^{-1}(1-u) \,du.
\end{eqnarray*}
The integrals exist because \mbox{$|G_N^{-1}(u)|, |F_N^{-1}(u)| \leq\max_{i=1}^N \max(|y_{1i}|, |y_{0i}|) < \infty$}.
\end{pf*}

Lemma~\ref{lemma:glivenko} below will be required in the proofs of
Propositions \ref{prop:asymptotics} and \ref{prop:rate}. In the special
case where the units of $U_N$ are independent and identically
distributed samples from a superpopulation, the first part of the lemma
reduces to the classical Glivenko--Cantelli theorem, and the
convergence implied by the second part follows from the conditional
bootstrap convergence results in van der Vaart and
Wellner \cite{VW1996}, Example~3.6.14. We
thank an anonymous reviewer for suggesting a more elegant way for
bounding (\ref{eq:Gvarbound}) and (\ref{eq:Fvarbound}) than our
original approach.

\begin{lemma}\label{lemma:glivenko}
Suppose conditions 1--3 of
Proposition~\ref{prop:asymptotics} hold. Then
\[
\sup_y \bigl|G(y)-G_N(y)\bigr| \rightarrow0\quad \mbox{and}\quad
\sup_y \bigl|F(y)-F_N(y)\bigr| \rightarrow0.
\]
In addition, given $\eta_1,\eta_0>0$, there exist two positive integers
$K_1(\eta_1)$ and $K_0(\eta_0)$ such that
\begin{eqnarray*}
\limsup_N \Bigl\{ N\pr\Bigl(\sup_y
\bigl|G(y)-\hat G_N(y)\bigr| \geq\eta_1\Bigr) \Bigr\}& \leq&
\frac{(1-\theta\rho)K_1(\eta_1)}{\theta\rho\eta_1^2},
\\
\limsup_N \Bigl\{ N\pr\Bigl(\sup_y
\bigl|F(y)-\hat F_N(y)\bigr| \geq\eta_0\Bigr) \Bigr\}& \leq&
\frac{ \{ 1 - \theta(1-\rho) \} K_0(\eta_0)}{\theta(1-\rho
)\eta
_0^2}.
\end{eqnarray*}
The integers are nonincreasing in $\eta$, and depend also on the
limiting distribution $H$ of $(y_1,y_0)$.
\end{lemma}

\begin{pf}
For the first part of the lemma, we follow the argument used in the
Glivenko--Cantelli theorem. Given $\eta_1>0$, there exists a partition
$-\infty= s_0 < s_1 < \cdots< s_{K_1(\eta_1)} = \infty$ such that
$G(s_i-) < G(s_{i-1}) + \eta_1/2$. For any $1 \leq i \leq K_1(\eta_1)$
and $s_{i-1} \leq s < s_i$,
\[
G(s_i-) - G_N(s_i-) - \eta_1/2
< G(s) - G_N(s) < G(s_{i-1}) - G_N(s_{i-1})
+ \eta_1/2,
\]
hence $\sup_y |G(y) - G_N(y)| < \eta_1$ if $|G(s_{i-1}) - G_N(s_{i-1})|
< \eta_1/2$ and $| G(s_i-) - G_N(s_i-) | < \eta_1/2$ for all $i$. By
conditions 2 and 3, this is
satisfied for all $N$ sufficiently large. The uniform convergence of
$F_N$ follows in the same way.

To establish the second part of the lemma, note that $\sup_y |G(y) -
\hat G_N(y) | < \eta_1$ on the set
\[
\bigcap_{i=1}^{K_1(\eta_1)} \bigl\{
\bigl|G(s_{i-1}) - \hat G_N(s_{i-1})\bigr|, \bigl|
G(s_i-) - \hat G_N(s_i-) \bigr| <
\eta_1/2 \bigr\}.
\]
Since $\pr\{ (\bigcap_i A_i)^c \} = \pr(\bigcup_i A_i^c) \leq\sum_i \pr
(A_i^c)$, we have
%
\begin{eqnarray}\label{eq:supCDFpbound}
&& \pr \Bigl\{ \sup_y \bigl|G(y)-\hat G_N(y)\bigr| \geq
\eta_1 \Bigr\}
\nonumber
\\
&&\qquad\leq \sum_{i=1}^{K_1(\eta_1)} \pr \bigl
\{\bigl|G(s_{i-1}) - \hat G_N(s_{i-1})\bigr| \geq
\eta_1/2 \bigr\}
\nonumber
\\
&&\qquad\quad{}+  \sum_{i=1}^{K_1(\eta_1)} \pr \bigl\{ \bigl|
G(s_i-) - \hat G_N(s_i-) \bigr| \geq
\eta_1/2 \bigr\}
\nonumber
\\[-8pt]
\\[-8pt]
\nonumber
&&\qquad\leq \sum_{i=1}^{K_1(\eta_1)} \pr\bigl\{\bigl |\hat
G_N(s_{i-1}) - G_N(s_{i-1})\bigr| \geq
\eta_1/2 - o(1) \bigr\}
\\
&&\qquad\quad{}+  \sum_{i=1}^{K_1(\eta_1)} \pr\bigl\{\bigl | \hat
G_N(s_i-) - G_N(s_i-)\bigr | \geq
\eta_1/2 - o(1) \bigr\}
\nonumber
\\
&&\qquad\leq \sum_{i=1}^{K_1(\eta_1)} \frac{ \Var_X \{ \hat G_N(s_{i-1})
\} +
\Var_X \{ \hat G_N(s_i-) \} }{ \{ \eta_1/2 - o(1) \}^2 },\nonumber
\end{eqnarray}
where the second inequality follows from $|G(y) - G_N(y)| = o(1)$, and
the last inequality from Chebyshev's inequality and the fact that $\E_X
\hat G_N(y) = G_N(y)$.

The argument used to derive (\ref{eq:varcalc}) can also be used to
bound the variances in (\ref{eq:supCDFpbound}). Noting that $\sigma
_N^2(I\{ y_1 \leq y \}) = G_N(y)\{ 1-G_N(y) \} \leq1/4$ and similarly
$\sigma_N^2(I\{ y_0 \leq y \}) \leq1/4$, we have for all $y$,
%
\begin{eqnarray}
\Var_X \hat G_N(y) &=& \frac{N-m}{(N-1)m}
\sigma_N^2\bigl(I\{ y_1 \leq y \}\bigr) =
\frac{N-m}{4(N-1)m}, \label{eq:Gvarbound}
\\
\Var_X \hat F_N(y) &=& \frac{N-(n-m)}{(N-1)(n-m)}
\sigma_N^2\bigl(I\{ y_0 \leq y \}\bigr) \leq
\frac{N-(n-m)}{4(N-1)(n-m)}. \label{eq:Fvarbound}
\end{eqnarray}
Plugging (\ref{eq:Gvarbound}) into (\ref{eq:supCDFpbound}) and taking
limits yields the desired result for $\hat G_N$, after absorbing the
factor of 2 into $K_1(\eta_1)$. The result for $\hat F_N$ can be
obtained in the same manner.
\end{pf}

\begin{pf*}{Proof of Proposition~\ref{prop:asymptotics}}
As indicated
in the proof outline, we proceed in several stages.

(i) \textit{Functional convergence of random distribution functions.}
Let $D([-\infty,\break \infty], \mathbb{R})^2$ be the Cartesian product of the
space of \textit{c\`adl\`ag} functions with itself, endowed with the
uniform metric induced by the norm $\| (v,u) \| = \max\{ \sup_y | v(y)
|,  \sup_y | u(y) | \}$. Thus, $D([-\infty,\infty],\mathbb{R})^2$
is a
nonseparable metric space. Lemma~\ref{lemma:glivenko} shows that the
distribution functions $(\hat G_N, \hat F_N)$ converge in probability
to $(G,F)$ in $D([-\infty,\infty],\mathbb{R})^2$. That is, $\pr(\|
(\hat G_N-G,\hat F_N-F) \| \geq\varepsilon) \rightarrow0$ for every
$\varepsilon>0$. As is the case with the lemma, the statement does not
require the use of outer measures because for each $N$, $(\hat G_N,
\hat F_N)$ can take on at most ${{N}\choose{n}} {{n}\choose{m}}$
distinct values in $D([-\infty,\infty],\mathbb{R})^2$; therefore, $\|
(\hat G_N-G,\hat F_N-F) \|$ is finite discrete valued.

(ii) \textit{Existence of random distributions $(\hat G_N',\hat F_N')$
defined on a common probability space.} Since the limit $(G,F)$ is
deterministic, the support of the limiting probability measure on
$D([-\infty,\infty],\mathbb{R})^2$ is a singleton. Applying the
Skorohod representation \cite{VW1996},\vspace*{1pt} Theorem~1.10.3, to $(\hat
G_N,\hat F_N)$ yields new random elements $(\hat G_N',\hat F_N')$ on
$D([-\infty,\infty],\mathbb{R})^2$ that have the same law as $(\hat
G_N,\hat F_N)$. Furthermore, $(\hat G_N',\hat F_N')$ converges to
$(G,F)$ almost everywhere, in the sense that along each sample path
$\omega'$ (in a set of measure one), the distribution functions
converge uniformly:
\[
\sup_y \bigl| \hat G_N'\bigl(y;
\omega'\bigr)-G(y) \bigr| \rightarrow0\quad \mbox{and}\quad \sup
_y\bigl | \hat F_N'\bigl(y;
\omega'\bigr)-F(y) \bigr| \rightarrow0.
\]

(iii) \textit{Convergence of $\E_{\hat G_N} (y_1^p)$ and $\E_{\hat F_N}
(y_0^p)$ for \textit{p} = 1, 2.} Define $\E_Q$ as the expectation
operator with respect to a distribution $Q$. Under condition 1, there exists $N_0$ such that $1/m \leq2/(\theta\rho
N)$ and $1/(n-m) \leq2/\{\theta(1-\rho) N\}$ for $N\geq N_0$. Then for
each $N \geq N_0$ and every realization of $(\hat G_N, \hat F_N)$,
condition 4 implies that as $\beta\rightarrow\infty$,
\begin{eqnarray*}
\E_{\hat G_N} \bigl( y_1^2I\bigl
\{y_1^2\geq\beta\bigr\} \bigr)& =& \sum
_{i:y_{1i}^2\geq
\beta
}^N \frac{ X^T_i y_{1i}^2 }{m} \leq\frac{2}{\theta\rho}
\sup_{N
\geq
N_0} \Biggl\{ \sum_{i:y_{1i}^2\geq\beta}^N
\frac{y_{1i}^2}{N} \Biggr\} \rightarrow0,
\\
\E_{\hat F_N} \bigl( y_0^2I\bigl
\{y_0^2\geq\beta\bigr\} \bigr) &=& \sum
_{i:y_{0i}^2\geq
\beta}^N \frac{ X^C_i y_{0i}^2 }{n-m} \leq\frac{2}{\theta(1-\rho)}
\sup_{N \geq N_0} \Biggl\{ \sum_{i:y_{0i}^2\geq\beta}^N
\frac
{y_{0i}^2}{N} \Biggr\} \rightarrow0.
\end{eqnarray*}
Recall that both ($\hat G_N',\hat F_N'$) and ($\hat G_N,\hat F_N$)
share the same finite discrete distribution. Thus, for almost all
sample paths $\omega'$ in the probability space of $(\hat G_N',\hat
F_N')$, the sequences of distributions represented by $\{ \hat
G_N'(\cdot;\omega') \}_N$ and $\{ \hat F_N'(\cdot;\omega') \}_N$ are
uniformly square-integrable. Moreover, since $\hat G_N'(\cdot;\omega')
\rightarrow G(\cdot)$, the random moments $\{ \E_{\hat G_N'} (y_1),
\E
_{\hat G_N'} (y_1^2) \}$ converge to $\{ \E_H (y_1),\break \E_H (y_1^2) \}$
almost everywhere, with the limits being finite. Similarly, $\{ \E
_{\hat F_N'} (y_0),\break \E_{\hat F_N'} (y_0^2) \} \rightarrow\{ \E_H
(y_0),  \E_H (y_0^2) \}$ almost everywhere as well. Translating this
back into convergence in probability for the first two random moments
of $\hat G_N$ and $\hat F_N$, we have
%
\begin{eqnarray}
\hat\S_N^2(y_1) & =& \frac{N-1}{N}
\frac{m}{m-1} \bigl[ \E_{\hat G_N} \bigl(y_1^2
\bigr) - \bigl\{\E_{\hat G_N} (y_1) \bigr\}^2 \bigr]
\rightarrow \Var_H(y_1), \label{eq:y1varplim}
\\
\hat\S_N^2(y_0) & = &\frac{N-1}{N}
\frac{n-m}{n-m-1} \bigl[ \E _{\hat
F_N} \bigl(y_0^2
\bigr) - \bigl\{\E_{\hat F_N} (y_0) \bigr\}^2 \bigr]
\rightarrow\Var_H(y_0) \label{eq:y0varplim}
\end{eqnarray}
in probability.

(iv) \textit{Convergence of $\hat\S_N^H(y_1,y_0)$ and $\hat\S
_N^L(y_1,y_0)$.} Define the distributions
$H^H(y_1,y_0) = \min\{ G(y_1),F(y_0) \}$ and $H^L(y_1,y_0) = \max\{
0,G(y_1)+F(y_0)-1 \}$, both of which have marginals $G$ and $F$. Using
Hoeffding's result from the proof of Lemma~\ref{lemma:covbound}, we
have that
\begin{eqnarray*}
\E_{H^H} (y_1y_0) &=& \sup_{h\in\mathcal{H}}
\E_h(y_1y_0),
\\
\E_{H^L} (y_1y_0) &=& \inf_{h\in\mathcal{H}}
\E_h(y_1y_0).
\end{eqnarray*}
Now fix a sample path and define two sequences of distributions $\hat
H_N^{H'}(y_1,y_0;\break \omega') = \min\{ \hat G_N'(y_1;\omega'),\hat
F_N'(y_0;\omega') \}$ and $\hat H_N^{L'}(y_1,y_0;\omega') = \max\{
0,\hat G_N'(y_1;\omega')+\hat F_N'(y_0; \omega')-1 \}$. It is clear that
$\hat H_N^{H'}(\cdot,\cdot; \omega')$ converges to $H^H(\cdot,\cdot)$
and $\hat H_N^{L'}(\cdot,\cdot;\break\omega')$ converges to $H^{L}(\cdot
,\cdot
)$ pointwise. Given that the product $y_1y_0$ is also uniformly
integrable with respect to almost all sequences $\{ \hat H_N^{H'}(\cdot
,\cdot;\omega') \}_N$ and\vadjust{\goodbreak} $\{ \hat H_N^{L'}(\cdot,\cdot;\omega') \}_N$
because $\{|XY|\geq\beta^2 \} \subset\{|X|\geq\beta\} \cup\{
|Y|\geq\beta\}$, it follows that $\E_{\hat H_N^{H'}} (y_1 y_0)
\rightarrow\sup_{h\in\mathcal{H}} \E_h(y_1y_0)$ and $\E_{\hat
H_N^{L'}} (y_1 y_0) \rightarrow\inf_{h\in\mathcal{H}} \E_h(y_1y_0)$
almost everywhere. Thus,
%
\begin{eqnarray}
\hat\S_N^H(y_1,y_0)& =&
\E_{\hat H_N^H} (y_1 y_0) - \E_{\hat G_N}
(y_1) \E_{\hat F_N} (y_0) \rightarrow\sup
_{h\in\mathcal{H}} \Cov _h(y_1,y_0),
\label{eq:covUBplim}
\\
\hat\S_N^L (y_1,y_0) &=&
\E_{\hat H_N^L} (y_1 y_0) - \E_{\hat G_N}
(y_1) \E_{\hat F_N} (y_0) \rightarrow\inf
_{h\in\mathcal{H}} \Cov_h(y_1,y_0)
\label{eq:covLBplim}
\end{eqnarray}
in probability. Plugging (\ref{eq:y1varplim})--(\ref{eq:covLBplim})
into (\ref{eq:varbound}) then yields the proposition.
\end{pf*}

\begin{proposition}\label{prop:rate}
Suppose conditions 1--3 of
Proposition~\ref{prop:asymptotics} hold, and that $y_1$ and $y_0$ are
bounded: $|y_{1i}|,|y_{0i}| \leq C < \infty$ for all $i$. Given
$\varepsilon>0$, for any $\varepsilon_1,\varepsilon_2,\varepsilon_3>0$
such that $\sum_i \varepsilon_i = \varepsilon$,
\begin{eqnarray*}
\pr\bigl( N\bigl| \hat V^H_N - V^H_N
\bigr|  \geq\varepsilon\bigr) &\leq&\mathcal{O} \biggl( \frac{C^4}{N}
\kappa_1(\varepsilon_1) \biggr),
\\
\pr\bigl( N\bigl| \hat V^L_N - V^L_N
\bigr|  \geq\varepsilon\bigr) &\leq&\mathcal{O} \biggl( \frac{C^4}{N} \bigl\{ 1/
\varepsilon_1^2 + \kappa_2(
\varepsilon_2) + \kappa _3(\varepsilon_3)
\bigr\} \biggr),
\end{eqnarray*}
where $\kappa_1(\varepsilon_1)$, $\kappa_2(\varepsilon_2)$ and
$\kappa
_3(\varepsilon_3)$ depend on the limiting distribution $H$.
\end{proposition}

\begin{pf} Define the bivariate distribution functions
$H^H_N(y_1,y_0) =\break  \min\{ G_N(y_1), F_N(y_0) \}$, $H^L_N(y_1,y_0) \vadjust{\goodbreak}= \max
( 0,G_N(y_1)+F_N(y_0)-1 )$,\break  $\hat H^H_N(y_1, y_0) = \min\{ \hat
G_N(y_1), \hat F_N(y_0) \}$, and $\hat H^L_N(y_1,y_0) = \max( 0, \hat
G_N(y_1)+\break  \hat F_N(y_0)-1 )$. Let $\E_Q$ be the expectation operator
with respect to a distribution $Q$. Using another result by Hoeffding
as recounted in Lehmann \cite{lehmann}, Lemma~2, the following covariances can
be expressed as
\begin{eqnarray*}
\hat\S_N^H(y_1,y_0) & = &\int
_{-C}^C\int_{-C}^C
\hat H_N^H(y_1,y_0)
\,dy_1 \,dy_0 \\
&&{}- \int_{-C}^C
\hat G_N(y_1) \,dy_1 \int_{-C}^C
\hat F_N(y_0) \,dy_0
\\
& =&\int_{-C}^C\int_{-C}^C
\hat H_N^H(y_1,y_0)
\,dy_1 \,dy_0 - C^2 + C \E _{\hat G_N}
(y_1)
\\
&&{} + C \E_{\hat F_N} (y_0) - \E_{\hat G_N} (y_1)
\E_{\hat F_N} (y_0),
\\
\S_N^H(y_1,y_0) & = &\int
_{-C}^C \int_{-C}^C
H_N^H(y_1,y_0) \,dy_1
\,dy_0 \\
&&{}- \int_{-C}^C
G_N(y_1) \,dy_1 \int_{-C}^C
F_N(y_0) \,dy_0
\\
& =& \int_{-C}^C \int_{-C}^C
H_N^H(y_1,y_0) \,dy_1
\,dy_0 - C^2 + C \E_{G_N} (y_1)
\\
&&{} + C \E_{F_N} (y_0) - \E_{G_N} (y_1)
\E_{F_N} (y_0),
\end{eqnarray*}
where the second equality follows from the identity $\E(W) = C - \int_{-C}^C \pr(W\leq w)\,dw$ for any random variable $W$ bounded by $C$. Then
\begin{eqnarray*}
N\bigl(\hat V_N^H - V_N^H\bigr)
& = &\frac{N-m}{m-1} \biggl\{ \E_{\hat G_N} \bigl(y_1^2
\bigr) - \frac{N(m-1)}{m(N-1)} \E_{G_N} \bigl(y_1^2
\bigr) \biggr\}
\\
&&{} - \frac{N-m}{m-1} \biggl[ \bigl\{ \E_{\hat G_N} (y_1) \bigr
\}^2 - \frac
{N(m-1)}{m(N-1)} \bigl\{ \E_{G_N} (y_1)
\bigr\}^2 \biggr]
\\
&&{} + \frac{N-(n-m)}{n-m-1} \biggl\{ \E_{\hat F_N} \bigl(y_0^2
\bigr) - \frac
{N(n-m-1)}{(N-1)(n-m)} \E_{F_N} \bigl(y_0^2
\bigr) \biggr\}
\\
& &{}- \frac{N-(n-m)}{n-m-1} \biggl[ \bigl\{ \E_{\hat F_N} (y_0)
\bigr\}^2 - \frac
{N(n-m-1)}{(N-1)(n-m)} \bigl\{ \E_{F_N}
(y_0) \bigr\}^2 \biggr]
\\
&&{} + \frac{2N}{N-1} C \bigl\{ \E_{\hat G_N} (y_1) -
\E_{G_N} (y_1) + \E _{\hat F_N}(y_0) -
\E_{F_N} (y_0) \bigr\}
\\
& &{}- \frac{2N}{N-1} \bigl\{ \E_{\hat G_N} (y_1)
\E_{\hat F_N}(y_0) - \E _{G_N} (y_1)
\E_{F_N} (y_0) \bigr\}
\\
&&{} + \frac{2N}{N-1} \int_{[-C,C]^2} \bigl\{ \hat
H_N^H(y_1,y_0) -
H_N^H(y_1,y_0) \bigr\}
\,dy_1 \,dy_0.
\end{eqnarray*}
To obtain the desired result, we proceed by bounding the probability
that each of the seven terms are large. Let $\nu_1,\ldots,\nu_8 > 0$ be
a tuple whose sum is $\varepsilon_1$. For the first term,
\begin{eqnarray*}
&& \pr \biggl\{ \frac{N-m}{m-1} \biggl\llvert \E_{\hat G_N}
\bigl(y_1^2\bigr) - \frac
{N(m-1)}{m(N-1)} \E_{G_N}
\bigl(y_1^2\bigr) \biggr\rrvert \geq\nu_1
\biggr\}
\\
&&\qquad\leq \pr \biggl\{ \bigl\llvert \E_{\hat G_N} \bigl(y_1^2
\bigr) - \E_{G_N} \bigl(y_1^2\bigr) \bigr\rrvert
\geq\frac{(m-1)\nu_1}{N-m} - \frac{(N-m)C^2}{m(N-1)} \biggr\}
\\
&&\qquad\leq \Var_X \bigl\{\E_{\hat G_N} \bigl(y_1^2
\bigr) \bigr\}\Big/ \biggl\{ \frac
{(m-1)\nu_1}{N-m} - o(1) \biggr\}^2
\\
&&\qquad\leq \frac{(N-m) C^4}{(N-1)m} \Big/ \biggl\{ \frac{(m-1)\nu_1}{N-m} - o(1) \biggr
\}^2,
\end{eqnarray*}
where the first inequality follows from $\llvert  \E_{\hat G_N} (y_1^2) -
\beta_N \E_{G_N} (y_1^2) \rrvert  \leq\llvert  \E_{\hat G_N} (y_1^2)
- \E
_{G_N} (y_1^2) \rrvert  + \llvert  (1 - \beta_N) \rrvert  \E_{G_N}
(y_1^2)$, and the second inequality from Chebyshev's inequality and the
fact that $\E_X \E_{\hat G_N} (y_1^p) = \E_{G_N} (y^p)$. The bound on
the variance is obtained in the same way as (\ref{eq:Gvarbound}). Thus,
%
\begin{eqnarray}
\label{eq:pbound1} &&\limsup_N N \pr \biggl\{ \frac{N-m}{m-1}
\biggl\llvert \E_{\hat G_N} \bigl(y_1^2\bigr) -
\frac{N(m-1)}{m(N-1)} \E_{G_N} \bigl(y_1^2\bigr)
\biggr\rrvert \geq\nu_1 \biggr\}
\nonumber
\\[-8pt]
\\[-8pt]
\nonumber
&&\qquad\leq \frac{(1-\theta\rho)^3 C^4}{\theta^3 \rho^3 \nu_1^2}.
\end{eqnarray}
For the second term,
\begin{eqnarray*}
&& \pr \biggl\{ \frac{N-m}{m-1} \biggl\llvert \bigl\{ \E_{\hat G_N}
(y_1) \bigr\}^2 - \frac{N(m-1)}{m(N-1)} \bigl\{
\E_{G_N} (y_1) \bigr\}^2 \biggr\rrvert \geq
\nu_2 \biggr\}
\\
&&\qquad\leq \pr \biggl\{ \biggl\llvert \E_{\hat G_N} (y_1) - \biggl\{
\frac
{N(m-1)}{m(N-1)} \biggr\}^{1/2} \E_{G_N} (y_1)
\biggr\rrvert \\
&&\hspace*{55pt}{}\geq\frac
{(m-1)\nu_2/  \{ C(N-m)  \} }{ 1 +  \{
{N(m-1)}/{(m(N-1))}  \}^{1/2} } \biggr\}
\\
&&\qquad\leq \pr \biggl\{ \bigl\llvert \E_{\hat G_N} (y_1) -
\E_{G_N} (y_1) \bigr\rrvert \geq\frac{(m-1)\nu_2}{ \{ 2+o(1) \}
 (N-m) C} - o(1)
\biggr\}
\\
&&\qquad\leq \Var_X \bigl\{ \E_{\hat G_N} (y_1) \bigr\} \Big/
\biggl[ \frac
{(m-1)\nu_2 }{ \{2+o(1)\} (N-m) C} - o(1) \biggr]^2
\\
&&\qquad \leq\frac{(N-m) C^2}{(N-1)m} \Big/ \biggl[ \frac{(m-1)\nu_2 }{ \{
2+o(1)\}
(N-m) C} - o(1)
\biggr]^2,
\end{eqnarray*}
where the first inequality follows from the identity $u^2-v^2 =
(u+v)(u-v)$. Hence,
%
\begin{eqnarray}
\label{eq:pbound2} &&\limsup_N N \pr \biggl\{ \frac{N-m}{m-1}
\biggl\llvert \bigl\{ \E_{\hat G_N} (y_1) \bigr\}^2 -
\frac{N(m-1)}{m(N-1)} \bigl\{ \E_{G_N} (y_1) \bigr
\}^2 \biggr\rrvert \geq \nu_2 \biggr\}
\nonumber
\\[-8pt]
\\[-8pt]
\nonumber
&&\qquad\leq \frac{4(1-\theta\rho)^3 C^4 }{\theta^3\rho^3 \nu_2^2}.
\end{eqnarray}
The same arguments apply to the third and fourth terms:
%
\begin{eqnarray}\label{eq:pbound3}
 &&\limsup_N N \pr \biggl\{ \frac{N-(n-m)}{n-m-1} \biggl\llvert
\E_{\hat F_N} \bigl(y_0^2\bigr) - \frac{N(n-m-1)}{(N-1)(n-m)}
\E_{F_N} \bigl(y_0^2\bigr) \biggr\rrvert \geq
\nu _3 \biggr\}
\nonumber\hspace*{-25pt}
\\[-8pt]
\\[-8pt]
\nonumber
&&\quad\leq \frac{ \{ 1-\theta(1-\rho) \}^3 C^4}{\theta^3 (1-\rho)^3
\nu_3^2},\hspace*{-25pt}\\
 \label{eq:pbound4}
 &&\limsup_N N \pr\biggl\{ \frac{N-(n-m)}{n-m-1} \biggl\llvert
\bigl\{ \E_{\hat
F_N} (y_0) \bigr\}^2 -
\frac{N(n-m-1)}{(N-1)(n-m)} \bigl\{ \E_{F_N} (y_0) \bigr
\}^2 \biggr\rrvert \geq\nu_4 \biggr\}
\nonumber\hspace*{-25pt}
\\[-8pt]
\\[-8pt]
\nonumber
&&\quad\leq \frac{4 \{ 1-\theta(1-\rho) \}^3 C^4 }{\theta^3(1-\rho)^3
\nu_4^2}.\hspace*{-25pt}
\end{eqnarray}
For the fifth term,
\begin{eqnarray*}
&& \pr \biggl[ \frac{2NC}{N-1} \bigl\llvert \E_{\hat G_N}
(y_1) - \E_{G_N} (y_1) + \E_{\hat F_N}(y_0)
- \E_{F_N} (y_0) \bigr\rrvert < \nu_5 +
\nu_6 \biggr]
\\
&&\qquad\geq \pr \biggl\{ \bigl\llvert \E_{\hat G_N} (y_1) -
\E_{G_N} (y_1) \bigr\rrvert < \frac{(N-1)\nu_5}{2NC}, \\
&&\hspace*{50pt}\bigl
\llvert \E_{\hat F_N}( y_0) - \E_{F_N}
(y_0) \bigr\rrvert < \frac{(N-1)\nu_6}{2NC} \biggr\}
\\
&&\qquad\geq 1 - \pr \biggl\{ \bigl\llvert \E_{\hat G_N} (y_1) -
\E_{G_N} (y_1) \bigr\rrvert \geq\frac{(N-1)\nu_5}{2NC} \biggr\}
\\
&&\qquad\quad - \pr \biggl\{ \bigl\llvert \E_{\hat F_N} (y_0) -
\E_{F_N} (y_0) \bigr\rrvert \geq\frac{(N-1)\nu_6}{2NC} \biggr\}
\\
&&\qquad\geq 1 - \Var_X \bigl\{ \E_{\hat G_N} (y_1) \bigr
\}\Big/ \biggl\{ \frac
{(N-1)\nu_5}{2NC} \biggr\}^2 - \Var_X \bigl
\{ \E_{\hat F_N} (y_0) \bigr\} \Big/ \biggl\{ \frac{(N-1)\nu_6}{2NC}
\biggr\}^2,
\end{eqnarray*}
so we have
%
\begin{eqnarray}\quad
\label{eq:pbound5}&& \limsup_N N \pr \biggl[ \frac{2NC}{N-1}
\bigl\llvert \E_{\hat G_N} (y_1) - \E _{G_N}
(y_1) + \E_{\hat F_N}(y_0) - \E_{F_N}
(y_0) \bigr\rrvert \geq\nu _5 + \nu_6
\biggr]
\nonumber
\\[-8pt]
\\[-8pt]
\nonumber
&&\qquad\leq \frac{4(1-\theta\rho) C^4}{\theta\rho\nu_5^2} + \frac{4 \{
1-\theta(1-\rho) \} C^4}{\theta(1-\rho)\nu_6^2}.
\end{eqnarray}

For the sixth term, we use the fact that $|uv-u'v'| \leq|uv-u'v| +
|u'v-u'v'|$ to obtain
\begin{eqnarray*}
&& \pr \biggl\{ \frac{2N}{N-1} \bigl\llvert \E_{\hat G_N}(
y_1) \E_{\hat
F_N}(y_0) - \E_{G_N}
(y_1) \E_{F_N} (y_0) \bigr\rrvert <
\nu_7 + \nu_8 \biggr\}
\\
&&\qquad\geq 1 - \pr \biggl\{ \bigl\llvert \E_{\hat G_N} (y_1) -
\E_{G_N} (y_1) \bigr\rrvert \geq\frac{(N-1)\nu_7}{2NC} \biggr\}
\\
&&\qquad\quad{} - \pr \biggl\{ \bigl\llvert \E_{\hat F_N} (y_0) -
\E_{F_N} (y_0) \bigr\rrvert \geq\frac{(N-1)\nu_8}{2NC} \biggr\}.
\end{eqnarray*}
Following the rest of the derivation of (\ref{eq:pbound5}) gives
%
\begin{eqnarray}
\label{eq:pbound6}&& \limsup_N N \pr \biggl\{ \frac{2N}{N-1}
\bigl\llvert \E_{\hat G_N} (y_1) \E _{\hat F_N}(y_0)
- \E_{G_N} (y_1) \E_{F_N} (y_0) \bigr
\rrvert \geq\nu_7 + \nu_8 \biggr\}
\nonumber
\\[-8pt]
\\[-8pt]
\nonumber
&&\qquad\leq \frac{4(1-\theta\rho) C^4}{\theta\rho\nu_7^2} + \frac{4 \{
1-\theta(1-\rho) \} C^4 }{\theta(1-\rho)\nu_8^2}.
\end{eqnarray}

To bound the probability that the last term exceeds $1-\varepsilon_1$,
first note that $| u - u' | < \eta$ and $| v - v' | < \eta$ implies $|
\min(u,v) - \min(u',v') | < \eta$. This gives the third inequality below:
\begin{eqnarray*}
&& \pr \biggl[ \frac{2N}{N-1} \biggl\llvert \int_{[-C,C]^2}
\bigl\{ \hat H_N^H(y_1,y_0) -
H_N^H(y_1,y_0) \bigr\}
\,dy_1 \,dy_0 \biggr\rrvert < 1-\varepsilon_1
\biggr]
\\
&&\qquad\geq \pr \biggl\{ \sup_{y_1,y_0} \bigl\llvert \hat
H_N^H(y_1,y_0) -
H_N^H(y_1,y_0) \bigr\rrvert <
\frac{(N-1)(1-\varepsilon_1)}{8NC^2} \biggr\}
\\
&&\qquad\geq \pr \biggl[ \sup_{y_1,y_0} \bigl\llvert \min\bigl\{\hat
G_N(y_1), \hat F_N(y_0) \bigr\}
- \min\bigl\{G_N(y_1), F_N(y_0)
\bigr\} \bigr\rrvert \\
&&\hspace*{194pt}{}< \frac
{(N-1)(1-\varepsilon_1)}{8NC^2} \biggr]
\\
&&\qquad\geq \pr \biggl\{ \sup_y\bigl | G_N(y) - \hat
G_N(y) \bigr|, \sup_y \bigl| F_N(y) - \hat
F_N(y) \bigr| < \frac{(N-1)(1-\varepsilon_1)}{8NC^2} \biggr\}
\\
&&\qquad\geq \pr \biggl\{ \sup_y \bigl| G(y) - \hat G_N(y)
\bigr|, \sup_y \bigl| F(y) - \hat F_N(y)\bigr | <
\frac{(N-1)(1-\varepsilon_1)}{8NC^2} - o(1) \biggr\}
\\
&&\qquad\geq 1 - \pr \biggl\{ \sup_y\bigl | G(y) - \hat
G_N(y) \bigr| \geq\frac
{(N-1)(1-\varepsilon_1)}{8NC^2} - o(1) \biggr\}
\\
&&\qquad\quad{}-  \pr \biggl\{ \sup_y \bigl| F(y) - \hat F_N(y) \bigr|
\geq\frac
{(N-1)(1-\varepsilon_1)}{8NC^2} - o(1) \biggr\}.
\end{eqnarray*}
The fourth inequality follows from Lemma~\ref{lemma:glivenko} which
shows that $\sup_y |G(y)-G_N(y)| = o(1)$ and $\sup_y |F(y)-F_N(y)| =
o(1)$. We can now apply the second part of Lemma~\ref{lemma:glivenko}
to bound the probability above. Given $\xi> 0$,
\begin{eqnarray*}
&& \limsup_N N \pr \biggl[ \frac{2N}{N-1} \biggl\llvert
\int_{[-C,C]^2} \bigl\{ \hat H_N^H(y_1,y_0)
- H_N^H(y_1,y_0) \bigr\}
\,dy_1 \,dy_0 \biggr\rrvert \geq 1-\varepsilon_1
\biggr]
\\
&&\qquad\leq \limsup_N N \pr \biggl\{ \sup_y
\bigl| G(y) - \hat G_N(y) \bigr| \geq \frac
{1-\varepsilon_1}{8C^2} - \xi \biggr\}
\\
&&\qquad\quad{}+  \limsup_N N \pr \biggl\{ \sup_y \bigl|
F(y) - \hat F_N(y) \bigr| \geq\frac
{1-\varepsilon_1}{8C^2} - \xi \biggr\}
\\
&&\qquad\leq \frac{(1-\theta\rho) K_1  (
 ({(1-\varepsilon_1)}/{(8C^2)}) -
\xi ) }{ \theta\rho \{ ({(1-\varepsilon_1)}/{(8C^2)}) -
\xi
\}^2 } \\
&&\qquad\quad{}+ \frac{ \{ 1-\theta(1-\rho) \} K_0
 ( ({(1-\varepsilon_1)}/{(8C^2)}) - \xi ) }{ \theta(1-\rho)  \{
({(1-\varepsilon_1)}/{(8C^2)}) - \xi \}^2 }.
\end{eqnarray*}\eject\noindent
Since $\xi$ is arbitrary and both $K_1(\cdot)$ and $K_0(\cdot)$ are
nonincreasing, there exists $\kappa_1(\varepsilon_1)$ such that
%
\begin{eqnarray}
\label{eq:pbound7} \quad&& \limsup_N N \pr \biggl[ \frac{2N}{N-1}
\biggl\llvert \int_{[-C,C]^2} \bigl\{ \hat H_N^H(y_1,y_0)
- H_N^H(y_1,y_0) \bigr\}
\,dy_1 \,dy_0 \biggr\rrvert \geq 1-\varepsilon_1
\biggr]
\nonumber
\\[-8pt]
\\[-8pt]
\nonumber
&&\qquad\leq C^4 \kappa_1(\varepsilon_1).
\end{eqnarray}

The bounds (\ref{eq:pbound1})--(\ref{eq:pbound7}) imply that
\[
\limsup_N N\pr\bigl(N\bigl|\hat V_N^H-V_N\bigr|
\geq\varepsilon\bigr) \leq C^4 \Biggl( \sum
_{i=1}^8 \frac{c_i}{\nu_i^2} + \kappa_1(
\varepsilon_1) \Biggr).
\]
By minimizing the right-hand side over $\nu_1,\ldots,\nu_8 > 0$ subject
to the constraint $\nu_1+\cdots+\nu_8= \varepsilon_1$, the sum in the
parenthesis can be absorbed into $\kappa_1(\varepsilon_1)$, yielding
the desired convergence rate for $N\hat V_N^H$. To get the rate for
$N\hat V_N^L$, we repeat the argument used to derive (\ref
{eq:pbound7}). First note that $| u - u' | < \eta$ and $| v - v' | <
\zeta$ implies $| \max(0,u+v-1) - \max(0,u'+v'-1) | < \eta+ \zeta$.
This gives the second inequality below:
\begin{eqnarray*}
&& \pr \biggl[ \frac{2N}{N-1} \biggl\llvert \int_{[-C,C]^2}
\bigl\{ \hat H_N^L(y_1,y_0) -
H_N^L(y_1,y_0) \bigr\}
\,dy_1 \,dy_0 \biggr\rrvert < \varepsilon_2 +
\varepsilon _3 \biggr]
\\
&&\qquad\geq \pr \biggl[ \sup_{y_1,y_0} \bigl\llvert \max \bigl\{ 0,\hat
G_N(y_1)+\hat F_N(y_0)-1 \bigr\}
\\
&&\hspace*{65pt}{}- \max \bigl\{ 0,G_N(y_1)+F_N(y_0)-1
\bigr\} \bigr\rrvert
\\
&&\hspace*{97pt}\qquad\quad  < \frac{(N-1)(\varepsilon_2 + \varepsilon_3)}{8NC^2} \biggr]
\\
&&\qquad\geq \pr \biggl\{ \sup_y \bigl| G(y) - \hat G_N(y)
\bigr| < \frac
{(N-1)\varepsilon_2}{8NC^2} - o(1),
\\
&&\hspace*{17pt}\qquad\quad  \sup_y \bigl| F(y) - \hat F_N(y) \bigr| <
\frac{(N-1)\varepsilon
_3}{8NC^2} - o(1) \biggr\}
\\[-1pt]
&&\qquad\geq 1 - \pr \biggl\{ \sup_y \bigl| G(y) - \hat
G_N(y) \bigr| \geq\frac
{(N-1)\varepsilon_2}{8NC^2} - o(1) \biggr\}
\\[-1pt]
&&\qquad\quad{}-  \pr \biggl\{ \sup_y \bigl| F(y) - \hat F_N(y)\bigr |
\geq\frac
{(N-1)\varepsilon_3}{8NC^2} - o(1) \biggr\}.\vspace*{-1pt}
\end{eqnarray*}
Thus there exist $\kappa_2(\varepsilon_2)$ and $\kappa_3(\varepsilon
_3)$ such that
%
\begin{eqnarray}
\label{eq:pbound8} && \limsup_N N \pr \biggl[ \frac{2N}{N-1}
\biggl\llvert \int_{[-C,C]^2} \bigl\{ \hat H_N^L(y_1,y_0)
- H_N^L(y_1,y_0) \bigr\}
\,dy_1 \,dy_0 \biggr\rrvert \geq\varepsilon_2 +
\varepsilon_3 \biggr]
\nonumber
\\[-10pt]
\\[-10pt]
\nonumber
&&\qquad\leq C^4 \bigl\{ \kappa_2(\varepsilon_2) +
\kappa_3(\varepsilon_3) \bigr\}.
\end{eqnarray}
\upqed\end{pf}

\section{R code for implementing estimator} \label{rcode}

Here, we present\vspace*{1pt} R code for the function \texttt{sharp.var}, which
outputs the bound estimates $\hat V^H_N$ (given input \texttt
{upper=TRUE}) and $\hat V^L_N$ (given input\break  \texttt{upper=FALSE}). The
other inputs are \texttt{yt} (the observed outcomes under treatment),
\texttt{yc} (the observed outcomes under control) and \texttt{N} (the
total number of units in the population).

\mbox{}

{\footnotesize{
\begin{verbatim}
sharp.var <- function(yt,yc,N=length(c(yt,yc)),upper=TRUE) {
 m <- length(yt)
 n <- m + length(yc)
 FPvar <- function(x,N) (N-1)/(N*(length(x)-1))
   * sum((x - mean(x))^2)
 yt <- sort(yt)
 if(upper == TRUE) yc <- sort(yc) else
   yc <- sort(yc,decreasing=TRUE)
 p_i <- unique(sort(c(seq(0,n-m,1)/(n-m),seq(0,m,1)/m))) -
 .Machine$double.eps^.5
 p_i[1] <- .Machine$double.eps^.5
 yti <- yt[ceiling(p_i*m)]
 yci <- yc[ceiling(p_i*(n-m))]
 p_i_minus <- c(NA,p_i[1: (length(p_i)-1)])
 return(((N-m)/m * FPvar(yt,N) + (N-(n-m))/(n-m) * FPvar(yc,N)
   + 2*sum(((p_i-p_i_minus)*yti*yci)[2:length(p_i)])
   - 2*mean(yt)*mean(yc))/(N-1))
}
\end{verbatim}
}}
%

\section{Illustrative upper bound improvements} \label{illustrative}

In Table~\ref{tab3}, we present illustrations of the improvements in the
variance upper bounds by varying the marginal distributions of
potential outcomes over the Beta distribution family: the control
potential outcomes are assumed to be distributed according to
$\operatorname
{Beta}(\alpha_0,\beta_0)$, and the treatment potential outcomes
according to $\operatorname{Beta}(\alpha_1,\beta_1)$. Strictly
speaking, since
finite populations cannot have continuous marginals, the Beta\vadjust{\goodbreak}
distributions represent approximations to plausible marginals when $N$
is large. We report the ratios ${ V^H_N}/{ V^{a}_N}$ and ${ V^H_N}/ {
V^{b+}_N}$ (the limits of ${ \hat V^H_N}/{ \hat V^{a}_N}$ and ${ \hat
V^H_N}/ { \hat V^{b+}_N}$) under different values of ($\alpha_0,\beta
_0,\alpha_1,\beta_1$) while holding $m = n/2$ and $n = N$ fixed.

\begin{table}
\caption{Illustrative upper bound ratios given Beta distributed
potential outcomes}\label{tab3}
\begin{tabular*}{\textwidth}{@{\extracolsep{\fill}}lllllcc@{}}
\hline
& \multicolumn{1}{c}{$\bolds{\alpha_0}$} &
\multicolumn{1}{c}{$\bolds{\beta_0}$} & \multicolumn{1}{c}{$\bolds{\alpha_1}$} &
\multicolumn{1}{c}{$\bolds{\beta_1}$} & \multicolumn{1}{c}{$\bolds{{ V^H_N}/{
V^{a}_N}}$} & \multicolumn{1}{c@{}}{$\bolds{{ V^H_N}/ { V^{b+}_N}}$} \\
\hline
\phantom{0}1 & 0.1 & 0.1 & 0.1 & 0.1 & 1.00 & 1.00 \\
\phantom{0}2 & 0.1 & 0.1 & 0.1 & 1 & 0.68 & 0.79 \\
\phantom{0}3 & 0.1 & 0.1 & 0.1 & 2 & 0.61 & 0.81 \\
\phantom{0}4 & 0.1 & 0.1 & 1 & 1 & 0.92 & 0.97 \\
\phantom{0}5 & 0.1 & 0.1 & 1 & 2 & 0.86 & 0.95 \\
\phantom{0}6 & 0.1 & 0.1 & 2 & 2 & 0.86 & 0.96 \\
\phantom{0}7 & 1 & 1 & 0.1 & 0.1 & 0.92 & 0.97 \\
\phantom{0}8 & 1 & 1 & 0.1 & 1 & 0.81 & 0.84 \\
\phantom{0}9 & 1 & 1 & 0.1 & 2 & 0.71 & 0.83 \\
10 & 1 & 1 & 1 & 1 & 1.00 & 1.00 \\
11 & 1 & 1 & 1 & 2 & 0.98 & 0.99 \\
12 & 1 & 1 & 2 & 2 & 0.98 & 1.00 \\
13 & 2 & 2 & 0.1 & 0.1 & 0.86 & 0.96 \\
14 & 2 & 2 & 0.1 & 1 & 0.85 & 0.85 \\
15 & 2 & 2 & 0.1 & 2 & 0.76 & 0.83 \\
16 & 2 & 2 & 1 & 1 & 0.98 & 1.00 \\
17 & 2 & 2 & 1 & 2 & 0.99 & 0.99 \\
18 & 2 & 2 & 2 & 2 & 1.00 & 1.00 \\
\hline
\end{tabular*}\vspace*{-9pt}
\end{table}

Table~\ref{tab3} presents 18 scenarios, wherein $(\alpha_0,\beta_0) \in\{
(0.1,0.1),(1,1),(2,2)\}$, and $\alpha_1, \beta_1 \in\{0.1,1,2\}$. The
results are identical for $\operatorname{Beta}(\alpha_1,\beta_1)$
and\break  $\operatorname
{Beta}(\beta_1,\alpha_1)$; thus, we omit redundant results. The ratios
were computed via numerical quadrature using the \texttt{NIntegrate}
command in Mathematica 7.0.1.0 under the default settings.

Our results illustrate that when the marginal distributions are
identical (i.e., cases 1, 10 and 18), all upper bounds are identical,
since the Cauchy--Schwarz and AM-GM inequalities hold exactly. However,
as the marginal distributions diverge in shape (e.g., cases 3, 9 and
15), our proposed upper bound $ V^H_N$ materially outperforms Neyman's
bounds $ V^a_N$ and $ V^{b+}_N$.
\end{appendix}

\section*{Acknowledgements}
The authors thank Allison Carnegie, Ed Kaplan, Winston Lin, Cyrus
Samii, Aad van der Vaart and the review team for helpful comments.




\printaddresses

\end{document}